\documentclass{article}

\usepackage[a4paper,hmargin=3cm,vmargin=3cm]{geometry}

\usepackage{amsmath}
\usepackage{amsfonts}
\usepackage{epsfig}
\usepackage{theorem}
\usepackage[T1,T2A]{fontenc}
\usepackage[utf8]{inputenc}
\usepackage{comment}

\newtheorem{theorem}{Theorem}[section]

\newtheorem{proposition}[theorem]{Proposition}

{\theorembodyfont{\rmfamily} % No need in ``\rm'' command after ``\begin''

\newtheorem{example}[theorem]{Example}

\newtheorem{definition}[theorem]{Definition}

}

\def\Q{{\mathbb Q}}

\newcommand\al{\alpha}

\newcommand\la{\lambda}
\newcommand\be{\beta}
\def\eps{\varepsilon}

\newcommand\eq{\,=\,}

\newcommand\lb{\linebreak[4]}
\newcommand\pb{\pagebreak[4]}

\def\msk{\medskip}

\title{Enumeration of Weighted Plane Trees}

\author{Alexander K. Zvonkin\thanks{LaBRI, UMR 5800, Universit\'e de Bordeaux, 
Talence, 33400, France. E-mail: {\tt zvonkin@labri.fr}}}

\date{\today}

\begin{document}

\maketitle

\begin{abstract}
In weighted trees, all edges are endowed with positive integral weight. We enumerate 
weighted bicolored plane trees according to their weight and number of edges.
\end{abstract}

\section{Preliminaries}

This note is not intended for a journal publication: it does not contain difficult 
results, and the proofs use only standard and well-known techniques. Furthermore, 
some of the results are already known. However, they have their place in the context of 
the study of weighted trees, see \cite{PakZvo.1-14}, \cite{PakZvo.2-14}.

\begin{definition}[Weighted tree]\label{def:weighted}
A {\em weighted bicolored plane tree}, or a {\em weighted tree}, or just a
{\em tree}\/ for short, is a bicolored plane tree whose edges are endowed 
with positive integral {\em weights}. The sum of the weights of the edges 
of a tree is called the {\em total weight}\/ of the tree. 

The {\em degree}\/ of a vertex is the sum of the weights of the edges incident to this 
vertex. Obviously, the sum of the degrees of black vertices, as well as the sum of the
degrees of white vertices, is equal to the total weight $n$ of the tree. Let the tree
have $p$ black vertices, of degrees $\al_1,\ldots,\al_p$, and $q$~white vertices,
of degrees $\be_1,\ldots,\be_q$, respectively. Then the pair of partitions $(\al,\be)$,
$\al,\be\vdash n$, is called {\em passport}\/ of the tree.

The {\em weight distribution}\/ of a weighted tree is a partition 
$\mu\vdash n$, $\mu=(\mu_1,\mu_2,\ldots,\mu_m)$ where $m=p+q-1$ is 
the number of edges, and $\mu_i$, $i=1,\ldots,m$ are the weights 
of the edges. Leaving aside the weights and considering only the
underlying plane tree, we speak of a {\em topological tree}, which is a bicolored
plane tree. Weighted trees whose weight distribution is $\mu=1^n$ will 
be called {\em ordinary trees}: they coincide with the corresponding topological trees.

%We call a {\em leaf}\/ a vertex which has only one edge incident to it,
%whatever is the weight of this edge. By abuse of language, we will
%also call a leaf this edge itself.
\end{definition}

The adjective {\em plane}\/ in the above definition means that our trees
are considered not as mere graphs but as plane maps. More precisely,
this means that the cyclic order of branches around each vertex of the 
tree is fixed, and changing this order will in general give a different 
tree. {\em All the trees considered in this paper will be endowed with 
the ``plane'' structure}\/; therefore, the adjective ``plane'' will often 
be omitted.

\begin{figure}[htbp]
\begin{center}
\epsfig{file=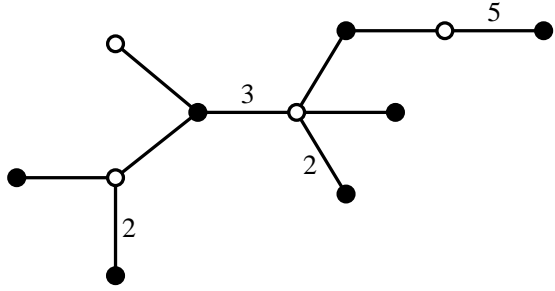,width=7.2cm}
\caption{\small Weighted bicolored plane tree. The weights which are not indicated
are equal to~1.}
\label{fig:wt}
\end{center}
\end{figure}

\begin{example}[Weighted tree]\label{ex:wt}
Figure \ref{fig:wt} shows an example of a weighted tree. The total weight of this tree 
is $n=18$; its passport is $(\al,\be)=(5^22^31^2,7^16^14^11^1)$; the weight 
distribution is \lb $\mu=5^13^12^21^6$.
\end{example}

\begin{definition}[Rooted tree]\label{def:rooted}
A tree with a distinguished edge is called a {\em rooted tree}, and the
distinguished edge itself is called its {\em root}. We consider the
root edge as being oriented from black to white.
\end{definition}

The goal of this paper is the enumeration of rooted weighted (bicolored plane) trees.

\section{Statement of the main theorem}

\begin{theorem}[Enumeration of weighted trees]\label{th:main}
Let $a_n$ be the number of rooted weighted \lb bicolored plane trees of
weight~$n$. Then the generating function $f(t)=\sum_{n\ge 0}a_nt^n$ is
equal to
\begin{eqnarray}
f(t) & = & \frac{1-t-\sqrt{1-6\,t+5\,t^2}}{2\,t} \nonumber \\
     & = & 1+t+3\,t^2+10\,t^3+36\,t^4+137\,t^5+543\,t^6+2219\,t^7+9285\,t^8
           +\ldots \label{gf:a_n}
\end{eqnarray}
Numbers $a_n$ satisfy the following recurrence relation:
\begin{eqnarray}\label{a_n-rec}
a_0=1, \quad a_1=1, \qquad a_{n+1}=a_n+\sum_{k=0}^n a_ka_{n-k} 
\quad \mbox{for} \quad n\ge 1.
\end{eqnarray}
The asymptotic formula for the numbers $a_n$ is
\begin{eqnarray}
a_n \sim \frac{1}{2}\sqrt{\frac{5}{\pi}}\cdot 5^n\,n^{-3/2}. \label{a_n-assymp}
\end{eqnarray}

Let $b_{m,n}$ be the number of rooted weighted bicolored plane trees
of weight~$n$ with $m$ edges. Then the generating function
$h(s,t)=\sum_{m,n\ge 0}b_{m,n}s^mt^n$ is equal to
\begin{eqnarray}
h(s,t) &=& \frac{1-t-\sqrt{1-(2+4s)\,t+(1+4s)\,t^2}}{2st} \nonumber \\
       &=& 1 + st + (s+2s^2)\,t^2 + (s+4s^2+5s^3)\,t^3 + 
             (s+6s^2+15s^3+14s^4)\,t^4 + \ldots \label{gf:b_m,n}
\end{eqnarray}
The following is an explicit formula for the numbers $b_{m,n}$:
\begin{eqnarray}
b_{m,n} \eq \binom{n-1}{m-1}\cdot{\rm Cat}_m 
       \eq \binom{n-1}{m-1}\cdot\frac{1}{m+1}\binom{2m}{m}, \label{b_m,n}
\end{eqnarray}
where ${\rm Cat}_m$ is the $m$th Catalan number.

Denote $|{\rm Aut}(T)|$ the order of the automorphism group of a tree~$T$. Let $c_n$ 
be the number of non-isomorphic non-rooted trees $T$ of weight $n$, each counted with the 
factor $1/|{\rm Aut}(T)|$. Then
\begin{eqnarray}
c_n \eq \sum_{T}\frac{1}{|{\rm Aut}(T)|} \eq \sum_{m=1}^n\frac{b_{m,n}}{m}\,, \label{c_n}
\end{eqnarray}
where the first sum is taken over all the non-isomorphic non-rooted trees $T$ of weight $n$.
\end{theorem}

The sequence $a_n$ is listed in the On-Line Encyclopedia of Integer Sequences \cite{OEIS} 
as the entry A002212. It has many interpretations, some of them coming from chemistry.
Among the various interpretations there are ``multi-trees'' (Roland Bacher, 2005) which 
correspond to our weighted trees. Almost all the above-stated formulas may also be found in
\cite{OEIS}.
%Thus, apparently, at least some part of the above results are 
%already known, but its seems that they remained non published\footnote{Roland Bacher 
%has told to the author that the enumeration of multi-trees is a byproduct of some 
%of his earlier topological studies.}.

\begin{figure}[htbp]
\begin{center}
\epsfig{file=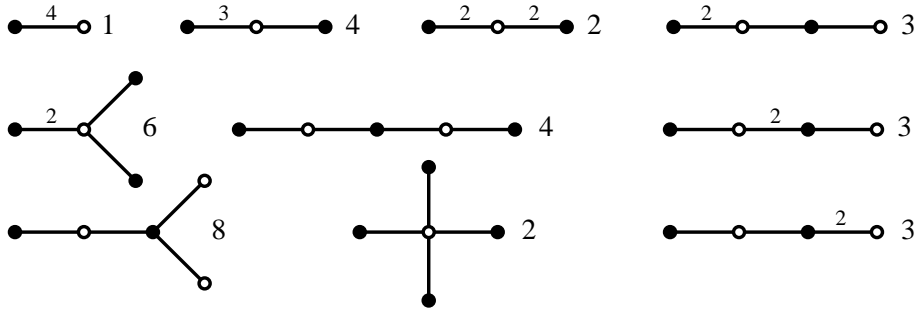,width=12cm}
\end{center}
\caption{\small Near each tree, the number of its possible rootings is indicated, 
with an eventual color exchange taken into account. The total number of rooted
trees is 36.}
\label{fig:ex-36}
\end{figure}

\begin{example}[Trees of weight 4]
Figure~\ref{fig:ex-36} shows the trees of weight~4. There are ten trees in the picture,
but in fact there are 16 non-isomorphic (non-rooted) trees of weight~4. Indeed, when we
exchange black and white, four trees remain isomorphic to themselves while six others 
don't, so we must add to the set the six missing trees. Near each tree, the number of 
its possible rootings is indicated, with color exchange taken into account. We see 
that the total number of trees is~36, which is the coefficient $a_4$ in front of $t^4$ 
in $f(t)$, see (\ref{gf:a_n}). Among these 36 trees, there is one tree with one edge, 
six trees with two edges, 15 trees with three edges, and 14 trees with four edges. 
These are the coefficients of the polynomial $s+6s^2+15s^3+14s^4$ which stands in front 
of $t^4$ in $h(s,t)$, see (\ref{gf:b_m,n}).

The number $c_4$, according to \eqref{c_n}, is equal to
$$
1 + \frac{6}{2} + \frac{15}{3} + \frac{14}{4} \eq 12\frac{1}{2}.
$$
And, indeed, among the 16 non-isomorphic non-rooted trees there are ten asymmetric 
trees, four trees with the symmetry order 2, and two trees with the symmetry order~4, 
which gives 
$$
10 + 4\cdot\frac{1}{2} + 2\cdot\frac{1}{4} \eq 12\frac{1}{2}.
$$
We leave the details to the reader.
\end{example}

\section{Dyck words and weighted Dyck words}

There is a standard way of encoding rooted topological (non-weighted) plane trees by 
{\em Dyck words}\/ and {\em Dyck paths}. We start on the left bank of the root edge 
and go around the tree in the clockwise direction, writing the letter $x$ when
we follow an edge for the first time, and the letter $y$ when we follow
it the second time on its opposite side. A Dyck path corresponding to a
Dyck word is a path on the plane which starts at the origin and takes
a step $(1,1)$ for every letter $x$ and a step $(1,-1)$ for every letter~$y$.
These objects may be easily characterized. 

For a word $w$ which is a concatenation of the three words, $w=u_1u_2u_3$ (either is 
allowed to be empty), we call $u_1$ a {\em prefix}, $u_2$ a {\em factor}, and $u_3$ 
a {\em suffix}\/ of~$w$.

\begin{definition}[Dyck words and Dyck paths]\label{def:dyck}
A {\em Dyck word}\/ is a word $w$ in the alphabet $\{x,y\}$ such that
$|w|_x=|w|_y$ (here $|w|_x$ and $|w|_y$ stand for the number of occurencies 
of $x$~and~$y$ in $w$), while for any prefix $u$ of $w$ we have 
$|u|_x\ge |u|_y$. A {\em Dyck path}\/ is a path on the plane which
starts at the origin, takes steps $(1,1)$ and $(1,-1)$, and finishes
on the horizontal axis, while always staying on the upper half-plane.
%The word and the path may be considered as {\em rooted}\/ at their first 
%letter or step. 
%A Dyck word is {\em prime}\/ if it is not a concatenation 
%of two or more non-empty Dyck words.
\end{definition}

\begin{figure}[htbp]
\begin{center}
\epsfig{file=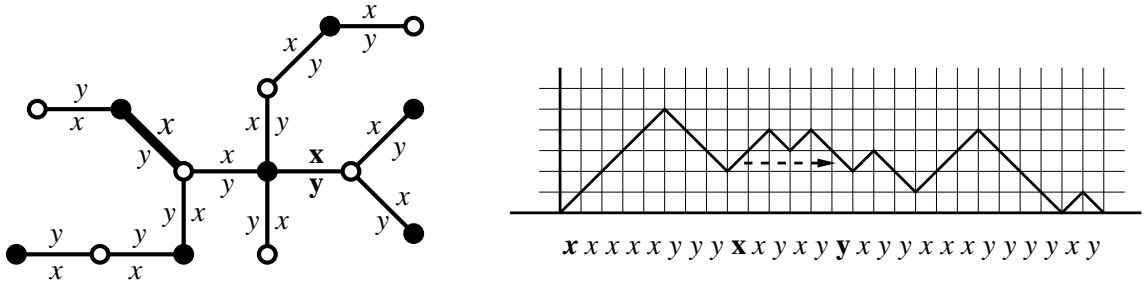,width=15cm}
\end{center}
\caption{\small A rooted bicolored plane tree and its encoding by the corresponding 
Dyck word and Dyck path.}
\label{fig:dyck}
\end{figure}

Figure~\ref{fig:dyck} illustrates these notions. The root edge in the tree
is shown by the thick line. Note that a Dyck word may be empty; then it corresponds
to the tree consisting of a single vertex. Making an exception to the general rule,
we do not color this vertex in black or white. Thus, there exists a single empty
word, and a single tree without edges.

\pb

The following proposition is a trivial consequence of the construction.

\begin{proposition}[Trees and Dyck words]\label{prop:tree-word}
There is a bijection between rooted bicolored plane trees and Dyck words.
\end{proposition}

There remains very little to do in order to describe weighted trees.
There is a natural notion of coupling of the letters of a Dyck word:
a {\em couple}\/ is the pair of letters $(x,y)$ standing on the 
opposite sides of the same edge. It is easy to recognize a couple in a word
or in the path. Let $x,y$ be a pair of letters in a Dyck word, where $x$ comes 
before $y$ in $w$. Consider the factor $xuy$ of $w$ which starts with $x$ and
terminates with $y$. Then the pair $(x,y)$ forms a couple if and only if the 
factor~$u$ between $x$ and $y$ is a Dyck word. In a Dyck path, we take an ascending
step corresponding to a letter $x$ and go horizontally until we meet a descending
step opposite to it: this step corresponds to the letter $y$ which forms a couple
with $x$. In Figure~\ref{fig:dyck}, an example of a couple, both on the tree and
in the word, is indicated in a boldface font, and the dashed arrow shows how to 
find the ``opposite'' step. The couple $(x,y)$ in which $x$ is the very first letter
of the Dyck word, corresponds to the root edge.

Now, returning to the weighted trees, we do the following: for every edge of the tree, 
we take the corresponding couple $(x,y)$ and replace it with $(x_i,y_i)$ where $i$ is
the weight of the edge.

\begin{definition}[Weighted Dyck words]\label{def:gen-dyck}
A {\em weighted Dyck word}\/ is a word in the infinite alphabet
$\{x_i,y_i\}_{i\ge 1}$ which is a Dyck word in which every couple 
of letters $(x,y)$ is replaced by a certain couple $(x_i,y_i)$.
We say that a couple $(x_i,y_i)$ has the weight $i$, and the weight
of a word is the sum of the weights of all its couples.
%All the notions of being rooted, prime, etc., remain unchanged.
\end{definition}

\begin{proposition}[Weighted trees and weighted Dyck words]
There is a bijection between rooted weighted bicolored plane trees and
weighted Dyck words.
\end{proposition}

\section{Proof of the main theorem}

Every non-empty Dyck word $w$ has a unique decomposition of the form $w=xuyv$ where 
$u$ and $v$ are themselves Dyck words (maybe empty). Here, obviously, $x$ is the first 
letter of~$w$, and $y$ is the letter coupled with it. The corresponding step in the 
Dyck path is the descending step of the first return of the path to the horizontal axis.

In the same way, every non-empty weighted Dyck word $w$ has a unique decomposition 
of the form $x_iuy_iv$ for some $i\ge 1$, where $u$ and $v$ are weighted Dyck words.

Let $\cal D$ be the formal sum of all the weighted Dyck words, that is, the formal 
power series  
\begin{eqnarray}\label{eq:series}
{\cal D} \eq \eps+x_1y_1+x_2y_2+x_1y_1x_1y_1+x_1x_1y_1y_1+x_3y_3+x_2y_2x_1y_1+
x_2x_1y_1y_2+\ldots
\end{eqnarray}
in non-commuting variables $x_i,y_i$, $i=1,2,\ldots$, where $\eps$ stands for the empty 
word. (In order to write down a series we must choose a total order on the words. 
A particular choice of the order is irrelevant. In \eqref{eq:series}, the words are
ordered according to their weight.) Then, the above decomposition of the words of 
$\cal D$ in the form $x_iuy_iv$ implies the following equation for $\cal D$:
\begin{eqnarray}\label{eq:eq-D}
{\cal D} \eq \eps+x_1{\cal D}y_1{\cal D}+x_2{\cal D}y_2{\cal D}+\ldots 
         \eq \eps+\sum_{i=1}^{\infty}x_i{\cal D}y_i{\cal D}\,.
\end{eqnarray}

%\pb

Now, do the following:
\begin{itemize}
\item   replace each letter $y_i$ in $\cal D$ by 1;
\item   replace each letter $x_i$ in $\cal D$ by $st^i$;
\item   make the variables $s$ and $t$ commute.
\end{itemize}

Then, every word $w$ in $\cal D$ is transformed into a word $s^mt^n$ where $m$ is 
the number of occurencies of the letters $x_i$, $i\ge 1$, in $w$ (or, equivalently, 
the number of edges of the weighted tree $T_w$ corresponding to $w$), and $n$ is 
the weight of $w$ (or, equivalently, the total weight of $T_w$). Therefore, combining 
similar terms we get the generating function $h(s,t)=\sum_{m,n\ge 0}b_{m,n}s^mt^t$.
At the same time, equation \eqref{eq:eq-D} is transformed into the following
quadratic equation for $h(s,t)$:
\begin{eqnarray}\label{eq:eq-h}
h \eq 1+s \left(\sum_{i=1}^{\infty}t^i\right) h^2 \eq 
1+\frac{st}{1-t}\cdot h^2.
\end{eqnarray}
Solving this equation, and choosing the sign in front of the square root
in such a way as to avoid a singularity at zero, we obtain 
formula~(\ref{gf:b_m,n}). Then, substituting $s=1$ in (\ref{gf:b_m,n})
we get (\ref{gf:a_n}). 

In order to obtain the asymptotic expression \eqref{a_n-assymp} for the numbers~$a_n$
it suffices to apply to~$f(t)$ the ready-made formulas of asymptotic analysis of the 
coefficients of generating functions -- see, for example, \cite{FlaSed-09}, Chapter VI. 
The only thing to note is that 
%the expression under the square root sign in $f$ is
$$
1-6t+5t^2 \eq (1-t)(1-5t).
$$ 

In order to prove \eqref{b_m,n}, we proceed as follows. There are ${\rm Cat}_m$ 
topological rooted trees with $m$ edges. Starting at the root edge, we go around 
a tree in the clockwise direction and attribute a non-zero weight to every newly
encountered edge. There are $\binom{n-1}{m-1}$ ways to do that. Indeed,
put $n$ dots in a row, and distribute $m-1$ separators among $n-1$ places between
the dots. This procedure splits the number $n$ into $m$ non-zero parts.

In order to prove the recurrence \eqref{a_n-rec}, consider separately the trees 
of weight $n+1$ having the root edge of weight~1, and the trees of weight $n+1$ 
having the root edge of weight $i\ge 2$. The weighted Dyck words corresponding 
to the trees of the first kind are of the form $x_1uy_1v$, where $u$~and~$v$ are 
themselves weighted Dyck words. The sum of the weights of $u$ and $v$ is $n$; 
denoting the weight of $u$ by $k$, so that the weight of $v$ becomes $n-k$, and 
summing over the $k=0,1,\ldots,n$, we get the term $\sum_{k=0}^n a_ka_{n-k}$ 
of~\eqref{a_n-rec}. Now, all the trees of weight $n+1$ having the root edge of
weight $i\ge 2$ are obtained from the trees of weight $n$ having the root edge 
of weight $i-1$, by adding one unit to the weight of the root. 
This gives the term $a_n$ in the right-hand part of \eqref{a_n-rec}.

Finally, \eqref{c_n} follows from the fact that there are $m$ choices of a root
edge in a tree $T$ with $m$ edges, but if this tree has non-trivial symmetries then
some of these choices produce isomorphic rooted trees. The number of non-isomorphic
rootings is $m/|{\rm Aut}(T)|$. Thus, dividing by $m$, we get the factor $1/|{\rm Aut}(T)|$.

\msk

Theorem \ref{th:main} is proved. \hfill$\Box$

\section{Enumeration of ordinary trees according to their passport}\label{sec:GJ}

Let $\la\vdash n$, $\la=(\la_1,\la_2,\ldots\la_k)$ be a partition of and integer $n$. 
Let us write $\la$ in the {\em power notation}\/:
$$
\la = 1^{d_1}2^{d_2}\ldots n^{d_n}, \qquad\mbox{where}\qquad
\sum_{i=1}^n d_i = k, \qquad \sum_{i=1}^n i\cdot d_i = n,
$$
so that $d_i$ is the number of the parts of $\la$ equal to $i$. Denote
\begin{eqnarray}\label{eq:N}
N(\la) \eq \frac{(k-1)!}{d_1!\,d_2!\,\ldots\, d_n!}\,.
\end{eqnarray}
The following theorem was proved in \cite{Tutte-64} (1964) and later generalized
in \cite{GouJac-92} (1992):

\begin{theorem}[Ordinary trees with a given passport]\label{th:GJ}
The number of rooted ordinary bicolored plane trees with the passport $(\al,\be)$ is
equal to 
\begin{eqnarray}\label{GJ-rooted}
nN(\al)N(\be). 
\end{eqnarray}
Respectively, the number of the non-isomorphic ordinary bicolored plane trees with the 
passport $(\al,\be)$, each one of them counted with the factor $1/|{\rm Aut}(T)|$, is
\begin{eqnarray}\label{GJ}
\sum_T \frac{1}{|{\rm Aut}\,(T)|} \eq N(\al)N(\be)
\end{eqnarray}
where the sum is taken over the the non-isomorphic ordinary bicolored plane trees 
with the passport $(\al,\be)$.
\end{theorem}

\begin{figure}[htbp]
\begin{center}
\epsfig{file=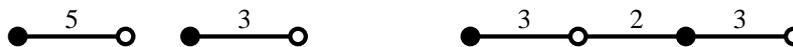,width=10.5cm}
\caption{\small The same passport $(5^13^1,5^13^1)$ is realized by a forest and by a tree.}
\label{fig:forest-and-tree}
\end{center}
\end{figure}

The above theorem is much more powerful than our Theorem \ref{th:main}. First of all,
it gives an explicit formula; and, what is more important, it enumerates the trees not
according to one or two parameters (as the weight and the number of edges in our case) 
but according to their passport. In the weighted case, a major difficulty in obtaining
a similar formula stems from the fact that the same passport can be realized by a tree 
and by a forest, as one can see in a very simple example of Figure~\ref{fig:forest-and-tree}.
Therefore, an inclusion-exclusion procedure might be unavoidable. This is indeed what 
takes place in~\cite{Kochetkov-13}.

\vspace{-5mm}

\begin{thebibliography}{99}
\bibitem{FlaSed-09}   {\bf Flajolet P., Sedgewick R.} 
         Analytic Combinatorics. --
         Cambridge Univ. Press, 2009, XIV+810~pp.
\bibitem{GouJac-92} {\bf Goulden I.\,P., Jackson D.\,M.}
        The combinatorial relationship between trees, cacti and certain
        connection coefficients for the symmetric group.
        {\em Europ. J. Combinat.}, 1992, vol.~{\bf 13}, 357--365.
\bibitem{Kochetkov-13} {\bf Kochetkov Yu.}
		Enumeration of one class of plane weighted trees. --
		\texttt{arXiv:1310.6208v1}, 11~pp.
\bibitem{OEIS} {\bf The On-Line Encyclopedia of Integer Sequences},
        \texttt{http://oeis.org/}.
\bibitem{PakZvo.1-14} {\bf Pakovich F., Zvonkin A.\,K.}
        Minimum degree of the difference of two polynomials over~$\Q$, 
        and weighted plane trees. -- \texttt{arXiv:1306.4141v1}. To appear in
		{\em Selecta Mathematica}, 2014.
\bibitem{PakZvo.2-14} {\bf Pakovich F., Zvonkin A.\,K.}
        Minimum degree of the difference of two polynomials over~$\Q$.
        Part~II: Davenport--Zannier triples. -- In preparation
        (a preliminary version may be found at
        \texttt{http://www.labri.fr/perso/zvonkin/}).
\bibitem{Tutte-64} {\bf Tutte W.\,T.}
		The number of planted plane trees with a given partition. --
		{\em Amer. Math. Monthly}, 1964, vol.~{\bf 71}, no.~3, 272--277.
\end{thebibliography}
\end{document}